\documentclass[12pt]{article}                                               

\input{amssym.def}                                                           
\input{amssym.tex}                                                           
                                                                             
\begin{document}                                                             
\title{Langlands reciprocity for  the even dimensional noncommutative tori}

\author{Igor  ~Nikolaev
\footnote{Partially supported 
by NSERC.}}


\date{}
 \maketitle


\newtheorem{thm}{Theorem}
\newtheorem{lem}{Lemma}
\newtheorem{dfn}{Definition}
\newtheorem{rmk}{Remark}
\newtheorem{cor}{Corollary}
\newtheorem{cnj}{Conjecture}
\newtheorem{exm}{Example}


\begin{abstract}
We conjecture an explicit  formula for the higher dimensional Dirichlet character; 
the formula is based on the $K$-theory of the so-called  noncommutative tori. 
It is proved,  that our conjecture is true for the two-dimensional
and one-dimensional (degenerate) noncommutative tori; in the second case,
one gets a noncommutative analog of the Artin reciprocity law.

\vspace{7mm}

{\it Key words and phrases:  Langlands program,  noncommutative tori}

\vspace{5mm}
{\it AMS  Subj. Class.:  11M55; 46L85}
\end{abstract}

\section*{Introduction}
The aim of the underlying note is to bring some evidence in favor of 
the following analog of the Langlands reciprocity \cite{Gel1}:
\begin{cnj}\label{cnj1}
{\bf (Langlands conjecture for noncommutative tori)}
Let $K$ be a finite extension of the  rational numbers ${\Bbb Q}$
with the  Galois group $Gal~(K | {\Bbb Q})$;
~for an irreducible representation $\sigma_{n+1}: Gal~(K | {\Bbb Q})\to GL_{n+1}({\Bbb C})$,
there exists a $2n$-dimensional  noncommutative torus with 
real  multiplication,    ${\cal A}_{RM}^{2n}$,  such that
$L(\sigma_{n+1},s)\equiv L({\cal A}_{RM}^{2n}, s)$, 
where $L(\sigma_{n+1}, s)$ is the Artin $L$-function and 
 $L({\cal A}_{RM}^{2n}, s)$ an $L$-function attached to the ${\cal A}_{RM}^{2n}$.
Moreover,  ${\cal A}_{RM}^{2n}$ is the image of
 an $n$-dimensional abelian variety $V_n(K)$ under
the  (generalized) Teichm\"uller functor $F_n$. 
\end{cnj}
For the notation and terminology we refer the reader to sections 1 and 3;  
the noncommutative torus ${\cal A}_{RM}^{2n}$  can  be regarded as 
a substitute  of  the ``automorphic cuspidal representation
$\pi_{\sigma_{n+1}}$ of the group $GL(n+1)$'' in terms of  the Langlands theory.
Roughly speaking, conjecture \ref{cnj1} says,  that the Galois extensions 
of the field of rational numbers come from the even dimensional noncommutative
tori with real multiplication.
Note, that the noncommutative  tori are intrinsic  to the problem, since they
classify the irreducible (infinite-dimensional) representations of 
the  Lie group $GL(n+1)$  \cite{Pog1};  such representations are the heart 
of the Langlands program \cite{Gel1}.  Our conjecture is supported by the following evidence. 
\begin{thm}\label{thm1}
Conjecture \ref{cnj1} is true for $n=1$ (resp., $n=0$) and $K$ abelian 
extension of an imaginary quadratic field $k$ (resp., the rational field ${\Bbb Q}$). 
\end{thm}
The structure of the note is as follows. A minimal necessary notation 
is introduced in section 1 and a brief summary of the Teichm\"uller 
functor(s) is given in section 3. Theorem \ref{thm1} is proved 
in section 2.

\section{Preliminaries}
\subsection{Noncommutative tori}
{\bf A. The $k$-dimensional noncommutative tori (\cite{Ell1},\cite{Rie1}).}
A {\it noncommutative $k$-torus} is the universal $C^*$-algebra
generated by $k$ unitary operators $u_1,\dots, u_k$; the operators
do not commute with each other, but their commutators 
$u_iu_ju_i^{-1}u_j^{-1}$ are fixed scalar multiples 
$\exp~(2\pi i\theta_{ij}), ~\theta_{ij}\in {\Bbb R}$ of the
identity operator.  The $k$-dimensional noncommutative torus, ${\cal A}_{\Theta}^k$,
is defined by a skew symmetric real matrix $\Theta=(\theta_{ij}),
~1\le i,j\le k$.  Further, we think of the ${\cal A}_{\Theta}^k$ as a noncommutative
topological space, whose algebraic $K$-theory yields  $K_0({\cal A}_{\Theta}^k)\cong {\Bbb Z}^{2^{k-1}}$ 
and  $K_1({\cal A}_{\Theta}^k)\cong {\Bbb Z}^{2^{k-1}}$.
The canonical trace $\tau$ on the $C^*$-algebra
${\cal A}_{\Theta}^k$ defines a homomorphism from 
$K_0({\cal A}_{\Theta}^k)$ to the real line ${\Bbb R}$;
under the homomorphism, the image of $K_0({\cal A}_{\Theta}^k)$
is a ${\Bbb Z}$-module, whose generators $\tau=(\tau_i)$ are polynomials 
in $\theta_{ij}$. (More precisely, $\tau=\exp (\Theta)$, where the
exterior algebra of $\theta_{ij}$ is nilpotent.) 
Recall, that the $C^*$-algebras ${\cal A}$ and ${\cal A}'$ are said
to be stably isomorphic (Morita equivalent), if ${\cal A}\otimes {\cal K}\cong
{\cal A}'\otimes {\cal K}$ for the $C^*$-algebra ${\cal K}$ of compact operators;
such an isomorphism indicates, that the $C^*$-algebras are homeomorphic as  noncommutative topological 
spaces.  By a result of Rieffel and Schwarz \cite{RiSch1}, the noncommutative
tori  ${\cal A}_{\Theta}^k$ and ${\cal A}_{\Theta'}^k$ are stably 
isomorphic,  if  the matrices $\Theta$ and $\Theta'$
belong to the same orbit of a subgroup $SO(k,k~|~{\Bbb Z})$ of the
group $GL_{2k}({\Bbb Z})$, which acts on $\Theta$ by the formula
$\Theta'=(A\Theta+B)~/~(C\Theta+D)$, where $(A, B,  C,  D)\in GL_{2k}({\Bbb Z})$
and  the matrices $A,B,C,D\in GL_k({\Bbb Z})$ satisfy the conditions:
\begin{equation}\label{eq1}
A^tD+C^tB=I,\quad A^tC+C^tA=0=B^tD+D^tB.
\end{equation}
(Here $I$ is the unit matrix and $t$ at the upper right of a matrix 
means a transpose of the matrix.)  
The group $SO(k, k ~| ~{\Bbb Z})$ can be equivalently defined as a
subgroup of the group  $SO(k, k ~| ~{\Bbb R})$ consisting of linear transformations 
of the space ${\Bbb R}^{2k}$,  which 
preserve the quadratic form $x_1x_{k+1}+x_2x_{k+2}+\dots+x_kx_{2k}$.

\medskip\noindent
{\bf B. The even dimensional normal tori.}
Further, we restrict to the case $k=2n$  (the even dimensional 
noncommutative tori).   It is known, that by the orthogonal linear transformations 
every (generic)  real even dimensional skew symmetric matrix can be brought to the  normal form:
\begin{equation}\label{eq2}
\Theta_0=
\left(
\matrix{      0 & \theta_1 &          &            &   \cr
      -\theta_1 &   0      &          &            &   \cr
                &          &  \ddots  &            &   \cr
                &          &          & 0          & \theta_n  \cr
                &          &          & -\theta_n  &    0
}\right)
\end{equation}
where $\theta_i>0$ are linearly independent over ${\Bbb Q}$. 
We shall consider the noncommutative tori ${\cal A}_{\Theta_0}^{2n}$,
given by the matrix (\ref{eq2}); we refer to the family 
as  a {\it normal family}. 
Recall, that any noncommutative torus has a canonical trace $\tau$,
which defines a homomorphism from $K_0({\cal A}_{\Theta}^k)\cong {\Bbb Z}^{2^{k-1}}$
to  ${\Bbb R}$; it follows from \cite{Ell1},  that the image of  
$K_0({\cal A}_{\Theta_0}^{2n})$ under the homomorphism has a
basis,  given by the formula
$\tau(K_0({\cal A}_{\Theta_0}^{2n}))=
{\Bbb Z}+\theta_1{\Bbb Z}+\dots+\theta_n{\Bbb Z}
+\sum_{i=n+1}^{2^{2n-1}}p_i(\theta){\Bbb Z}$,
where $p_i(\theta)\in {\Bbb Z}[1, \theta_1,\dots,\theta_n]$.

\medskip\noindent
{\bf C. The real multiplication (\cite{Man1}).}
The noncommutative torus ${\cal A}^k_{\Theta}$ is said to
have a {\it real multiplication}, if the endomorphism 
ring $End~(\tau(K_0({\cal A}_{\Theta}^k)))$ exceeds 
the trivial ring ${\Bbb Z}$.  Since any
endomorphism of the ${\Bbb Z}$-module  $\tau(K_0({\cal A}_{\Theta}^k))$
is the multiplication by a real number, it is easy to deduce, that
all the entries  of $\Theta=(\theta_{ij})$ are  algebraic integers.
(Indeed, the endomorphism is described by an integer matrix, which 
defines a polynomial equation involving $\theta_{ij}$.)
Thus, the noncommutative tori with real multiplication is 
a countable subset of all $k$-dimensional tori;  any element of the
set we shall denote by   ${\cal A}_{RM}^k$. 
Notice, that for the even dimensional normal tori with real multiplication,
the polynomials $p_i(\theta)$ produce the  algebraic 
integers in the extension of ${\Bbb Q}$ by $\theta_i$; any such
an integer is a linear combination (over ${\Bbb Z}$) of the $\theta_i$.
Thus,  the trace formula reduces to
$\tau(K_0({\cal A}_{RM}^{2n}))=
{\Bbb Z}+\theta_1{\Bbb Z}+\dots+\theta_n{\Bbb Z}$.

\subsection{$L$-function of noncommutative tori}
We consider even dimensional normal tori with real multiplication.
Denote by $A$ a positive integer matrix, whose (normalized) Perron-Frobenius 
eigenvector coincides with the vector $\theta=(1,\theta_1,\dots,\theta_n)$
and $A$ is not a power of a positive integer matrix;
in other words, $A\theta=\lambda_A\theta$, where $A\in GL_{n+1}({\Bbb Z})$
and $\lambda_A$ is the corresponding eigenvalue.  
(Explicitly,  $A$ can be obtained from vector $\theta$ as
the matrix of minimal  period of the  Jacobi-Perron continued fraction of $\theta$ 
\cite{B}.) Let $p$ be  a prime number;  take the matrix $A^p$ and consider its 
characteristic polynomial $char~(A^p)=x^{n+1}+a_1x^n+\dots+a_nx+1$. 
We introduce the following notation:
\begin{equation}\label{eq3}
L_p^{n+1}:=\left(
\matrix{
a_1     & a_2    & \dots  & a_n & p\cr
-1      & 0      & \dots  &  0      & 0\cr
\vdots  & \vdots  & \ddots & \vdots  & \vdots\cr
0       &  0     & \dots  &    -1   & 0
}\right). 
\end{equation}
A {\it local} zeta function of the ${\cal A}_{RM}^{2n}$
is defined as the reciprocal of  $det~(I_{n+1}-L_p^{n+1}z)$; 
in other words, 
\begin{equation}\label{eq4}
\zeta_p({\cal A}_{RM}^{2n}, z):= {1\over 1-a_1z+a_2z^2-\dots-a_n z^n +pz^{n+1}}, \quad z\in {\Bbb C}.
\end{equation}
An {\it $L$-function} of ${\cal A}_{RM}^{2n}$ is a product of the local
zetas over all,  but a finite number,  of primes
$L({\cal A}_{RM}^{2n},s)=\prod_{p ~\nmid  ~tr^2(A)-(n+1)^2}
\zeta_p({\cal A}_{RM}^{2n}, p^{-s}), ~s\in {\Bbb C}$.
\begin{rmk}
It will be shown, that for $n=0$ and $n=1$ formula (\ref{eq4}) fits conjecture \ref{cnj1};
for $n\ge 2$ it is an open problem based on an observation, that the crossed
product ${\cal A}_{RM}^{2n}\rtimes_{L_p^{n+1}}{\Bbb Z}$ is a proper noncommutative
analog of the (higher dimensional) Tate module,  where matrix $L_p^{n+1}$
corresponds to the Frobenius automorphism of the module \cite{S}, p.172.  
\end{rmk}

\section{Proof of theorem \ref{thm1}}
\subsection{Case $n=1$}
Each one-dimensional abelian variety is  a non-singular elliptic curve;
choose this curve to have  complex multiplication
by (an order in) the imaginary quadratic field $k$ and  denote
such a curve  by $E_{CM}$.  Then, by theory of complex multiplication,
the (maximal) abelian extension of $k$ coincides with 
the minimal field of definition of the curve $E_{CM}$, 
i.e. $E_{CM}\cong E(K)$ \cite{S}.
The Teichm\"uller functor $F:=F_1$  maps $E(K)$ into
a two-dimensional  noncommutative torus with real multiplication (section 3);
we shall denote the torus by ${\cal A}_{RM}^2$.  To  calculate the corresponding 
$L$-function $L({\cal A}_{RM}^2,s)$, let  $A$ be a $2\times 2$  positive integer matrix,
whose  normalized Perron-Frobenius eigenvenctor is  $(1, \theta_1)$. 
For a prime $p$,  the characteristic polynomial of the matrix $A^p$ writes as 
$char~(A^p)=x^2+tr~(A^p)x+1$ and the matrix $L_p^2$ takes the
form:
\begin{equation}\label{eq5}
L_p^2=\left(\matrix{tr~(A^p) & p\cr -1 & 0}\right).
\end{equation}
The corresponding local zeta function  
$\zeta_p({\cal A}_{RM}^2,z)=(1-tr~(A^p)z+pz^2)^{-1}$. 
We have to prove, that   $\zeta_p({\cal A}_{RM}^2,z)=\zeta_p(E_{CM},z)$, where 
$\zeta_p(E_{CM},z)$ is the local zeta function for
the elliptic curve $E_{CM}$;  the proof will be arranged into
a series of lemmas \ref{lm1}-\ref{lm5}.

Recall, that $\zeta_p(E_{CM},z)=(1-tr~(\psi_{E(K)}({\goth P}))z+pz^2)^{-1}$,
where $\psi_{E(K)}$ is the Gr\"ossencharacter on $K$, ${\goth P}$ the
prime ideal of $K$ over $p$ and  $tr$ is the trace of algebraic number
\cite{S}, Ch.2, \S 9.
Roughly, our proof consists in construction of representation $\rho$
of $\psi_{E(K)}$ into  the group of invertible elements (units) 
of $End~(\tau(K_0({\cal A}_{RM}^2)))$,  such  that  
$tr~(\psi_{E(K)}({\goth P}))=tr~(\rho(\psi_{E(K)}({\goth P})))=tr~(A^p)$.
This will be achieved with the help of an explicit formula 
for the Teichm\"uller functor $F$ (\cite{Nik1}, p.524):
\begin{equation}\label{eq6}
F: \left(\matrix{a & b\cr c & d}\right)\in End~(E_{CM})
\longmapsto
 \left(\matrix{a & b\cr -c & -d}\right)\in End~({\Bbb A}_{RM}).
\end{equation}
\begin{lem}\label{lm1}
Let $A=(a,b,c,d)$ be an integer matrix with $ad-bc\ne 0$
and $b=1$. Then $A$ is similar to the matrix
$(a+d, 1, c-ad, 0)$.
\end{lem}
{\it Proof.} Indeed, consider a matrix $(1,0,d,1)\in SL_2({\Bbb Z})$;
it is verified directly, that the matrix realizes the required similarity. 
$\square$
\begin{lem}\label{lm2}
The matrix $A=(a+d, 1, c-ad, 0)$ is similar to its
transpose $A^t=(a+d, c-ad, 1, 0)$. 
\end{lem}
{\it Proof.} We shall use the following criterion: the
(integer) matrices $A$ and $B$ are similar, if and only if 
the characteristic matrices $xI-A$ and $xI-B$ have the same Smith normal
form.  The calculation for the matrix $xI-A$ gives:
$$
\left(\matrix{x-a-d & -1\cr ad-c & x}\right)\sim
\left(\matrix{x-a-d & -1\cr  x^2-(a+d)x+ad-c & 0}\right)\sim
$$
$$
\sim \left(\matrix{1 & 0\cr 0 & x^2-(a+d)x+ad-c}\right),
$$
where $\sim$ are the elementary operations between the rows (columns)  
of the matrix.  Similarly, a calculation for the matrix $xI-A^t$
gives:
$$
\left(\matrix{x-a-d & ad-c\cr -1 & x}\right)\sim
\left(\matrix{x-a-d & x^2-(a+d)x+ad-c\cr -1& 0}\right)\sim
$$
$$
\sim\left(\matrix{1 & 0\cr 0 & x^2-(a+d)x+ad-c}\right).
$$
Thus, $(xI-A)\sim (xI-A^t)$ and lemma \ref{lm2} follows.
$\square$
\begin{cor}\label{cr1}
The matrices $(a, 1, c, d)$ and $(a+d, c-ad, 1, 0)$
are similar. 
\end{cor}
Let $E_{CM}$ be elliptic curve with  complex multiplication
by an order $R$ in the ring of integers of the imaginary quadratic
field $k$. Then ${\Bbb A}_{RM}=F(E_{CM})$ is a noncommutative torus
with real multiplication by the order ${\goth R}$ in the ring of
integers of a  real quadratic field ${\goth k}$ (section 3). 
Let $tr~(\alpha)=\alpha+\bar\alpha$ be the
trace function of a (quadratic) algebraic number field. 
\begin{lem}\label{lm3}
Each $\alpha\in R$ goes under $F$ into an $\omega\in {\goth R}$,
such that $tr~(\alpha)=tr~(\omega)$. 
\end{lem}
{\it Proof.} Recall that each $\alpha\in R$ can be written  in a matrix
form for a given base $\{\omega_1,\omega_2\}$ of the lattice
$\Lambda$. Namely, $\alpha\omega_1 = a\omega_1 +b\omega_2$ and 
$\alpha\omega_2 = c\omega_1 +d\omega_2$,
where $(a,b,c,d)$ is an integer matrix with $ad-bc\ne 0$. 
Note that $tr~(\alpha)=a+d$ and $b\tau^2+(a-d)\tau-c=0$, where
$\tau=\omega_2/\omega_1$. Since $\tau$ is an algebraic integer,
we conclude that $b=1$. 
In view of corollary \ref{cr1}, in a base $\{\omega_1',\omega_2'\}$,
the $\alpha$ has a matrix form $(a+d, c-ad, 1, 0)$. 
To calculate $\omega\in {\goth R}$ corresponding to $\alpha$,
we apply formula (\ref{eq6}), which gives us:
\begin{equation}\label{eq7}
F: \left(\matrix{a+d & c-ad\cr 1 & 0}\right)
\longmapsto
 \left(\matrix{a+d & c-ad\cr -1 & 0}\right). 
\end{equation}
In a given base $\{\lambda_1,\lambda_2\}$ of
the pseudo-lattice ${\Bbb Z}+{\Bbb Z}\theta$ one can write
$\omega\lambda_1 = (a+d)\lambda_1 +(c-ad)\lambda_2$ and 
$\omega\lambda_2 = -\lambda_1$.
It is an easy exercise to verify that $\omega$ is a real 
quadratic integer with $tr~(\omega)=a+d$; the latter coincides
with the $tr~(\alpha)$.
$\square$

\bigskip\noindent
Let $\omega\in {\goth R}$ be an endomorphism of the  pseudo-lattice
${\Bbb Z}+{\Bbb Z}\theta$ of degree $deg~(\omega):=\omega\bar\omega=n$.
The endomorphism maps  pseudo-lattice to a sub-lattice of index $n$.
Any such has a form ${\Bbb Z}+(n\theta){\Bbb Z}$ \cite{BS}, p.131. 
Let us calculate $\omega$ in a base $\{1,n\theta\}$, when $\omega$ is
given by the matrix $(a+d, c-ad, -1, 0)$. In this case $n=c-ad$ and $\omega$
induces an automorphism $\omega^*=(a+d, 1, -1, 0)$ of the sublattice  
${\Bbb Z}+(n\theta){\Bbb Z}$ according to the matrix equation:  
\begin{equation}\label{eq8}
\left(\matrix{a+d & n\cr -1 & 0}\right)
\left(\matrix{1\cr \theta}\right)=
\left(\matrix{a+d & 1\cr -1 & 0}\right)
\left(\matrix{1\cr n\theta}\right).
\end{equation}
Thus, one gets a map $\rho: {\goth R}\to {\goth R}^*$ given
by the formula $\omega=(a+d, n, -1, 0)\mapsto \omega^*=(a+d, 1, -1, 0)$,
where ${\goth R}^*$ is the group of units of  ${\goth R}$.
Since  $tr~(\omega^*)=a+d=tr~(\omega)$ and 
$\omega^*=\rho(\omega)$, one gets the following 
\begin{cor}\label{cr2}
For all $\omega\in {\goth R}$, it holds  $tr~(\omega)=tr~(\rho(\omega))$. 
\end{cor}
Note,  that ${\goth R}^*=\{\pm\varepsilon^k ~|~ k\in {\Bbb Z}\}$, where $\varepsilon>1$ is a 
fundamental unit of the order ${\goth R}\subseteq O_{\goth k}$; here  $O_{\goth k}$ means
the ring  of integers of a real quadratic  field ${\goth k}={\Bbb Q}(\theta)$. 
Choosing a sign in front of $\varepsilon^k$,  the following index map is defined
$\iota: R \buildrel\rm F\over\longrightarrow{\goth R}\buildrel\rm\rho \over
\longrightarrow {\goth R}^*\longrightarrow {\Bbb Z}$.
Let $\alpha\in R$ and $deg~(\alpha)=-n$. To calculate the
$\iota(\alpha)$,  recall some notation from Hasse
\cite{H}, \S 16.5.C. Let ${\Bbb Z}/n {\Bbb Z}$ be a cyclic group of
order $n$. For brevity, let $I={\Bbb Z}+{\Bbb Z}\theta$ be
a pseudo-lattice and $I_n={\Bbb Z}+(n\theta){\Bbb Z}$ its
sub-lattice of index $n$; the fundamental units of $I$
and $I_n$ are $\varepsilon$ and $\varepsilon_n$, respectively.
By ${\goth G}_n$ one understands a subgroup of ${\Bbb Z}/n {\Bbb Z}$
of prime residue classes $mod~n$. The ${\goth g}_n\subset {\goth G}_n$
is a subgroup of the non-zero divisors of the ${\goth G}_n$.
Finally, let $g_n$ be the smallest number, such that it divides 
$|{\goth G}_n/{\goth g}_n|$ and  $\varepsilon^{g_n}\in I_n$.
(The notation drastically simplifies in the case $n=p$ is a prime
number.)
\begin{lem}\label{lm4}
$\iota(\alpha)=g_n$. 
\end{lem}
{\it Proof.} Notice, that $deg~(\omega)=-deg~(\alpha)=n$,
where $\omega=F(\alpha)$. Then the map $\rho$ defines $I$ and $I_n$;
one can now apply the calculation of \cite{H}, pp 296-300.
Namely, Theorem $XIII^{\prime}$ on p. 298 yields the required
result. (We kept the notation of the original.) 
$\square$ 
\begin{cor}\label{cr3}
$\iota(\psi_{E(K)}({\goth P}))=p$. 
\end{cor}
{\it Proof.} It is known, that $deg~ (\psi_{E(K)}({\goth P}))=-p$,
where $\psi_{E(K)}({\goth P})\in R$ is the Gr\"ossencharacter.
To calculate the $g_n$ in the case $n=p$, notice that the 
${\goth G}_p\cong {\Bbb Z}/p{\Bbb Z}$ and ${\goth g}_p$ is trivial. 
Thus, $|{\goth G}_p/{\goth g}_p|=p$ is divisible only by $1$ or $p$.
Since $\varepsilon^1$ is not in $I_n$, one concludes that $g_p=p$.
The corollary follows.
$\square$ 
\begin{lem}\label{lm5}
$tr~(\psi_{E(K)}({\goth P}))=tr~(A^p)$. 
\end{lem}
{\it Proof.} It is not hard to see, that $A$ is a hyperbolic
matrix with the eigenvector $(1,\theta)$; the corresponding 
(Perron-Frobenius) eigenvalue is a fundamental unit $\varepsilon>1$ 
of the pseudo-lattice ${\Bbb Z}+{\Bbb Z}\theta$. In other words, $A$
is a matrix form of the algebraic number $\varepsilon$. It is immediate,
that $A^p$ is the matrix form for the $\varepsilon^p$ and $tr~(A^p)=tr~(\varepsilon^p)$.  
In view of lemma \ref{lm3} and corollary \ref{cr3},   
$tr~(\alpha)=tr~(F(\alpha))=tr~(\rho(F(\alpha)))$ for $\forall\alpha\in R$. 
In particular, if $\alpha=\psi_{E(K)}({\goth P})$ then, by corollary \ref{cr3}, 
one gets $\rho(F(\psi_{E(K)}({\goth P})))=\varepsilon^p$.  Taking traces in the 
last equation, we obtain the conclusion of lemma \ref{lm5}.
$\square$

\bigskip\noindent
The fact $\zeta_p({\cal A}_{RM}^2,z)=\zeta_p(E_{CM},z)$ follows from lemma \ref{lm5},
since the trace of the Gr\"ossencharacter coincides with such for the matrix $A^p$.
Let $E_{CM}$ be an elliptic curve with complex multiplication
by an order in the imaginary quadratic field $k$ and  $K$
 the minimal field of definition of the $E_{CM}$.
\begin{lem}\label{lm6}
$L(E_{CM},s)\equiv L(\sigma_2, s)$,
where $L(E_{CM},s)$ is the Hasse-Weil $L$-function of $E_{CM}$ and   
$L(\sigma_2, s)$ the Artin $L$-function for an irreducible  representation 
$\sigma_2: Gal~(K | k)\to GL_2({\Bbb C})$.
\end{lem}
{\it Proof.}  By the Deuring theorem (see e.g. \cite{S}, p.175),\linebreak
$L(E_{CM}, s)=L(\psi_K, s)L(\overline{\psi}_K, s)$, where
$L(\psi_K,s)$ is the Hecke $L$-series attached to the 
Gr\"ossencharacter $\psi: {\Bbb A}_K^*\to {\Bbb C}^*$; here 
${\Bbb A}_K^*$ denotes the adele ring of the field $K$ and the bar
means a complex conjugation. Notice, that since our elliptic curve has complex 
multiplication,  the group $Gal~(K | k)$ is abelian;  one can apply Theorem 5.1
\cite{Kna1}, which says that the Hecke $L$-series 
$L(\sigma_1\circ \theta_{K|k}, s)$ equals the Artin $L$-function $L(\sigma_1, s)$,
where $\psi_K=\sigma\circ \theta_{K|k}$ is the Gr\"ossencharacter and 
 $\theta_{K|k}: {\Bbb A}_K^*\to Gal~(K | k)$  the canonical homomorphism.
Thus, one gets $L(E_{CM}, s)\equiv L(\sigma_1,s)L(\overline{\sigma}_1,s)$,
where $\overline{\sigma}_1: Gal~(K | k)\to {\Bbb C}$ means a (complex) conjugate
representation of the Galois group. 
Consider the local factors of the Artin $L$-functions $L(\sigma_1,s)$ and $L(\overline{\sigma}_1,s)$;
it is immediate, that they are $(1-\sigma_1(Fr_p)p^{-s})^{-1}$ and  $(1-\overline{\sigma}_1(Fr_p)p^{-s})^{-1}$,
respectively. Let us consider a representation $\sigma_2: Gal~(K | k)\to GL_2({\Bbb C})$, 
such that 
\begin{equation}\label{eq9}
\sigma_2(Fr_p)=
\left(\matrix{\sigma_1(Fr_p) & 0\cr 0 & \overline{\sigma}_1(Fr_p)}\right). 
\end{equation}
It can be verified, that $det^{-1}(I_2-\sigma_2(Fr_p)p^{-s})=
 (1-\sigma_1(Fr_p)p^{-s})^{-1}(1-\overline{\sigma}_1(Fr_p)p^{-s})^{-1}$,
i.e. $L(\sigma_2,s)=L(\sigma_1,s)L(\overline{\sigma}_1,s)$. 
Lemma \ref{lm6} follows.
$\square$

\bigskip\noindent
By  lemma \ref{lm6}, we conclude, that 
$L({\cal A}_{RM}^2,s)\equiv L(\sigma_2,s)$ for an irreducible
representation $\sigma_2: Gal~(K|k)\to GL_2({\Bbb C})$.
It remains to notice that  $L(\sigma_2,s)= L(\sigma_2',s)$,
where $\sigma_2': Gal~(K|{\Bbb Q})\to GL_2({\Bbb C})$ \cite{Art1},
\S 3.  Case $n=1$ of  theorem \ref{thm1} follows.
$\square$

\subsection{Case $n=0$}
When $n=0$, one gets a one-dimensional (degenerate) noncommutative
torus; such an object, ${\cal A}_{\Bbb Q}$, can be obtained from the 
$2$-dimensional torus ${\cal A}^2_{\theta}$ by forcing $\theta=p/q\in {\Bbb Q}$ be a rational
number (hence our notation). 
One can always assume $\theta=0$ and, thus,   
$\tau(K_0({\cal A}_{\Bbb Q}))={\Bbb Z}$. 
To calculate matrix $L_p^1$, notice that the group of automorphisms of 
the ${\Bbb Z}$-module  $\tau(K_0({\cal A}_{\Bbb Q}))={\Bbb Z}$
is trivial, i.e. is a  multiplication by $\pm 1$; 
hence our  $1\times 1$  (real) matrix $A$ is either $1$ or $-1$.
Since $A$ must be positive, we conclude, that $A=1$. However,
$A=1$ is not a prime matrix, if one allows the complex entries;
indeed, for any $N>1$  matrix $A'=\zeta_N$ gives us $A=(A')^N$,  
where $\zeta_N=e^{2\pi i\over N}$ is the $N$-th root of unity.
 Therefore, $A=\zeta_N$ and  $L_p^1=tr~(A^p)=A^p=\zeta_N^p$.
A degenerate noncommutative torus, corresponding to the matrix $A=\zeta_N$,
we shall write as ${\cal A}_{\Bbb Q}^N$; in turn, such a torus is the
image (under the Teichm\"uller functor) of a zero-dimensional 
abelian variety,  which we denote by $V_0^N$. 
Suppose that $Gal~(K|{\Bbb Q})$ is abelian and   
let $\sigma: Gal~(K|{\Bbb Q})\to {\Bbb C}^{\times}$ be a
homomorphism. Then, by the Artin reciprocity \cite{Gel1}, there
exists an integer  $N_{\sigma}$ and a  Dirichlet character 
$\chi_{\sigma}: ({\Bbb Z}/N_{\sigma} {\Bbb Z})^{\times}\to {\Bbb C}^{\times}$,
such that $\sigma(Fr_p)=\chi_{\sigma}(p)$;
choose our zero-dimensional variety be $V_0^{N_{\sigma}}$.
In view of the notation, $L_p^1=\zeta_{N_{\sigma}}^p$; on the other hand,
it is verified directly, that $\zeta_{N_{\sigma}}^p=e^{{2\pi i\over N_{\sigma}}p}=
\chi_{\sigma}(p)$. Thus, $L_p^1=\chi_{\sigma}(p)$.
To obtain a local zeta function, we substitute $a_1=L_p^1$ into the
formula  (\ref{eq4}) and get
\begin{equation}\label{eq10}
\zeta_p({\cal A}_{\Bbb Q}^{N_{\sigma}},z)= {1\over 1-\chi_{\sigma}(p)z},
\end{equation}
where $\chi_{\sigma}(p)$ is the Dirichlet character. 
Therefore, $L({\cal A}^{N_{\sigma}}_{\Bbb Q}, s)\equiv L(s,\chi_{\sigma})$ 
is the Dirichlet $L$-series;  such a series, by construction,  coincides 
with the Artin $L$-series of the representation $\sigma: Gal~(K|{\Bbb Q})\to {\Bbb C}^{\times}$.
Case $n=0$ of theorem \ref{thm1} follows.
$\square$

\section{Teichm\"uller functors}
Denote by $\Lambda$  a lattice of rank $2n$; recall, that an $n$-dimensional (principally polarized) 
abelian  variety, $V_n$,  is  the complex torus ${\Bbb C}^n/\Lambda$,  which admits
an embedding into a projective space  \cite{M}. 
\subsection{Abelian varieties of dimension $n=1$}
{\bf A.  Basic example.} 
Let  $n=1$ and consider the complex torus  $V_1\cong {\Bbb C}/({\Bbb Z}+\tau {\Bbb Z})$;
it always embeds (via the Weierstrass $\wp$ function) into a 
projective space as a non-singular elliptic curve. 
Let  ${\Bbb H}=\{\tau=x+iy\in {\Bbb C} ~|~ y>0\}$ be the upper half-plane
and $\partial {\Bbb H}=\{\theta\in {\Bbb R} ~|~ y=0\}$ its (topological) 
boundary. We identify $V_1(\tau)$ with the points of ${\Bbb H}$ and 
${\cal A}^2_{\theta}$ with the points of $\partial {\Bbb H}$. 
Let us show, that the boundary is natural; the latter means, 
that the action of the modular group $SL_2({\Bbb Z})$ extends to
the boundary,  where it coincides with the  stable isomorphisms 
of tori.  Indeed,  conditions (\ref{eq1}) are equivalent to 
\begin{equation}\label{eq11}
A=\left(\matrix{a & 0\cr 0 & a}\right),\quad
B=\left(\matrix{0 & b\cr -b & 0}\right),\quad
C=\left(\matrix{0 & -c\cr c & 0}\right),\quad
D=\left(\matrix{d & 0\cr 0 & d}\right),
\end{equation}
where $ad-bc=1$, $a,b,c,d\in {\Bbb Z}$ and 
$\Theta'=(A\Theta+B)/(C\Theta+D)= 
(0,  {a\theta+b\over c\theta+d},  -{a\theta+b\over c\theta+d},  0)$.
Therefore,  $\theta'=(a\theta+b)(c\theta+d)^{-1}$ for a matrix
$(a, b, c, d)\in SL_2({\Bbb Z})$.
Thus, the action of $SL_2({\Bbb Z})$ extends to the boundary $\partial {\Bbb H}$,
where it induces  stable isomorphisms of the noncommutative
tori.

\medskip\noindent
{\bf B. The Teichm\"uller functor (\cite{Nik1}).}
There exists a continuous map $F_1: {\Bbb H}\to
\partial {\Bbb H}$, which sends isomorphic complex tori to
the stably isomorphic noncommutative tori. An exact result is
this.  Let $\phi$ be a closed form on the torus, whose trajectories
define a measured foliation; according to the Hubbard-Masur 
theorem (applied to the complex tori), this foliation 
corresponds to a point $\tau\in {\Bbb H}$. The map 
$F_1: {\Bbb H}\to\partial {\Bbb H}$
is defined by the formula $\tau\mapsto\theta=\int_{\gamma_2}\phi/\int_{\gamma_1}\phi$,
where $\gamma_1$ and $\gamma_2$ are generators of the first homology of the
torus.    The following is true: (i) ${\Bbb H}=\partial {\Bbb H}\times (0,\infty)$
is a trivial fiber bundle, whose projection map coincides with $F_1$;
(ii) $F_1$ is a functor, which sends isomorphic complex tori to
the stably isomorphic noncommutative tori.  We shall refer to $F_1$
as the {\it Teichm\"uller functor}.  Recall, that the complex 
torus ${\Bbb C}/({\Bbb Z}+\tau {\Bbb Z})$ is said to have  a complex multiplication,
if the endomorphism ring of the lattice $\Lambda={\Bbb Z}+\tau {\Bbb Z}$
exceeds the trivial ring ${\Bbb Z}$;  the complex multiplication happens
if and only if $\tau$ is an algebraic number in an  imaginary quadratic field.
The following is true: $F_1(V_1^{CM})={\cal A}^2_{RM}$, where $V_1^{CM}$
is a torus with   complex multiplication.

\subsection{Abelian varieties of dimension $n\ge 1$}
{\bf A.  The Siegel upper half-space (\cite{M}).}
The space  ${\Bbb H}_n:=\{\tau=(\tau_i)\in {\Bbb C}^{{n(n+1)\over 2}}~|~ Im~(\tau_i)>0\}$
of symmetric $n\times n$ matrices with complex entries is called a {\it Siegel upper half-space};
the points of ${\Bbb H}_n$ are one-to-one with  the $n$-dimensional principally polarized abelian varieties.
Let $Sp(2n, {\Bbb R})$ be the symplectic group;  it acts on ${\Bbb H}_n$
by the linear fractional transformations
$\tau\to \tau'=(a\tau+b)/(c\tau +d)$,
where $(a, b,  c, d)\in Sp(2n, {\Bbb R})$
and $a,b,c$ and $d$ are the  $n\times n$ matrices with real entries.
The abelian varieties $V_n$ and $V_n'$ are isomorphic, if and only
if,  $\tau$ and $\tau'$  belong to the same orbit of the group $Sp(2n, {\Bbb Z})$;
the action is discontinuous on ${\Bbb H}_n$ \cite{M}, Ch.2, \S 4. 
Denote by $\Sigma_{2n}$ a space of the $2n$-dimensional normal noncommutative
tori.  The following lemma is critical. 
\begin{lem}\label{lm7}
$Sp(2n, {\Bbb R})\subseteq O(n,n | {\Bbb R})$.
\end{lem}
{\it Proof.} 
(i) The group $O(n,n ~| ~{\Bbb R})$ can be defined as a subgroup 
of $GL_2({\Bbb R})$, which preserves the quadratic form
$f(x_1,\dots,x_{2n})=x_1x_{n+1}+x_2x_{n+2}+\dots+x_nx_{2n}$ \cite{RiSch1}. 
We shall denote $u_i=x_1$ for $1\le i\le n$
and $v_i=x_i$ for $n+1\le i\le 2n$. Consider the following  skew symmetric 
bilinear form $q(u,v)=u_1v_{n+1}+\dots+u_nv_{2n}-u_{n+1}v_1-\dots-u_{2n}v_n$,
where $u,v\in {\Bbb R}^{2n}$. It is known, that each linear
substitution $g\in Sp(2n, {\Bbb R})$ preserves the form $q(u,v)$.
Since $q(u,v)=f(x_1,\dots,x_{2n})- u_{n+1}v_1-\dots-u_{2n}v_n$,
one concludes that $g$ also preserves the form  $f(x_1,\dots,x_{2n})$,
i.e. $g\in O(n,n | {\Bbb R})$. It is easy to see, that the inclusion
is proper except the case $n=1$, i.e. when $Sp(2, {\Bbb R})\cong O(1,1 | {\Bbb R})\cong SL_2({\Bbb R})$.
Lemma follows.

\medskip
(ii)  We wish to give a second proof of this important fact, which is based on the 
explicit formulas for the block matrices $A,B,C$ and $D$. 
The fact that a symplectic linear transformation preserves
the  skew symmetric bilinear form $q(u,v)$ can be written 
in a matrix form:
\begin{equation}\label{eq12}
\left(\matrix{a & b\cr c & d}\right)^t
\left(\matrix{0 & I\cr -I & 0}\right)
\left(\matrix{a & b\cr c & d}\right)=
\left(\matrix{0 & I\cr -I & 0}\right),
\end{equation}
where $t$ is the transpose of a matrix. Performing the matrix
multiplication, one gets the following matrix identities
$a^td-c^tb=I, ~a^tc-c^ta=0=b^td-d^tb$.
Let us show, that these identities imply the Rieffel-Schwarz identities
(\ref{eq1}) imposed on the matrices $A, B, C$ and $D$. Indeed, in view
of the formulas (\ref{eq11}), the Rieffel-Schwarz identities can be written as:
\begin{equation}\label{eq13}
\left\{
\begin{array}{ccc}
\left(\matrix{a^t & 0\cr 0 & a^t}\right)
\left(\matrix{d & 0\cr 0 & d}\right) 
&+
\left(\matrix{0 & c^t\cr -c^t & 0}\right)
\left(\matrix{0 & b\cr -b & 0}\right) 
&=
\left(\matrix{I & 0\cr 0 & I}\right)
\nonumber\\
\left(\matrix{a^t & 0\cr 0 & a^t}\right)
\left(\matrix{0 & -c\cr c & 0}\right) 
&+
\left(\matrix{0 & c^t\cr -c^t & 0}\right)
\left(\matrix{a & 0\cr 0 & a}\right) 
&=
\left(\matrix{0 & 0\cr 0 & 0}\right)
\nonumber\\
\left(\matrix{0 & -b^t\cr b^t & 0}\right)
\left(\matrix{d & 0\cr 0 & d}\right) 
&+
\left(\matrix{d^t & 0\cr 0 & d^t}\right)
\left(\matrix{0 & b\cr -b & 0}\right) 
&=
\left(\matrix{0 & 0\cr 0 & 0}\right).
\end{array}
\right.
\end{equation}
A step by step matrix multiplication in (\ref{eq13}) shows that the identities 
$a^td-c^tb=I, ~a^tc-c^ta=0=b^td-d^tb$ imply the identities (\ref{eq13}).
(Beware: the operation is not commutative.)
Thus, any symplectic transformation satisfies the
Rieffel-Schwarz identities, i.e. belongs to the group $O(n,n | {\Bbb R})$.
Lemma \ref{lm7} follows.
$\square$

\medskip\noindent
{\bf B. The generalized Teichm\"uller functors.}
By lemma \ref{lm7},  the action of  $Sp(2n, {\Bbb Z})$ on the ${\Bbb H}_n$
extends to the $\Sigma_{2n}$,  where it acts by stable isomorphisms
of  the noncommutative tori;  thus,  $\Sigma_{2n}$ is  a  natural boundary 
of the Siegel upper half-space ${\Bbb H}_n$.  However,  unless $n=1$,  the
$\Sigma_{2n}$ is {\it not} a  topological boundary of  ${\Bbb H}_n$. 
Indeed,  $dim_{\Bbb R}({\Bbb H}_n)=n(n+1)$ and  $dim_{\Bbb R}(\partial {\Bbb H}_n)=n^2+n-1$,
while  $dim_{\Bbb R}(\Sigma_{2n})=n$.  Thus,  $\Sigma_{2n}$ is
 an $n$-dimensional subspace of the topological boundary of 
 ${\Bbb H}_n$;  this subspace is 
everywhere dense in $\partial {\Bbb H}_n$, since   
the $Sp(2n, {\Bbb Z})$-orbit of an element of $\Sigma_{2n}$ 
is everywhere dense in $\partial {\Bbb H}_n$ \cite{RiSch1}. 
A (conjectural)  continuous map $F_n: {\Bbb H}_n\to\Sigma_{2n}$,
we shall call a  {\it generalized Teichm\"uller functor}.
The $F_n$ has the following properties: (i) it sends each pair
of isomorphic abelian varieties to a pair of the stably isomorphic
even dimensional normal tori; (ii) the range of $F_n$ on
the  abelian varieties with complex multiplication
consists of the noncommutative tori with real multiplication.
As explained, such a functor has been constructed only in 
the case $n=1$;  the difficulties in higher dimensions are due to the 
lack (so far) of a proper Teichm\"uller theory for  the abelian varieties
of dimension $n\ge 2$.

\bigskip\noindent
{\sf Acknowledgment.} 
I am grateful to the referee for thoughtful comments.



\vskip1cm

\textsc{The Fields Institute for Mathematical Sciences, Toronto, ON, Canada,  
E-mail:} {\sf igor.v.nikolaev@gmail.com}

\smallskip
{\it Current address: 616-315 Holmwood Ave., Ottawa, ON, Canada, K1S 2R2}

\end{document}